\input ./style/arxiv-general.cfg
\documentclass[aop,MSNbibl,seceqn,dvips]{arximspdf}
\makeatletter
   \@ifpackageloaded{graphicx}{}{\usepackage{graphicx}}
\makeatother

\doi{10.1214/14-AOP912} 
\volume{43}
\issue{4}
\pubyear{2015}
\firstpage{1577}
\lastpage{1593}

\makeatletter

\newcommand{\iint}{\int\!\!\!\int}
\newcommand{\eqref}[1]{(\ref{#1})}
\newcommand {\IN}{\mathbb{N}}  
\newcommand {\IR}{\mathbb{R}}  
\newcommand {\ImF}{\mathcal{F}}  
\newcommand {\ImG}{\mathcal{G}}  
\newcommand {\ImJ}{\mathcal{J}}  
\newcommand {\ImC}{\mathcal{C}}  
\newcommand {\ImD}{\mathcal{D}}  
\newcommand {\ImR}{\mathcal{R}}  
\newcommand {\mbt}{\mathbf{t}}  
\newcommand {\mbs}{\mathbf{s}}  
\newtheorem{theorem}{Theorem}[section]

\newtheorem{lemma}[theorem]{Lemma}
\newproclaim{prop}[theorem]{Proposition}
\newproclaim{defi}[theorem]{Definition}
\newproclaim{remark}[theorem]{Remark}
\makeatother

\begin{document}
\begin{frontmatter}

\title{Multiple points of the Brownian sheet in~critical~dimensions}
\runtitle{Multiple points of the Brownian sheet}

\begin{aug}
\author[A]{\fnms{Robert C.} \snm{Dalang}\corref{}\ead[label=e1]{robert.dalang@epfl.ch}\thanksref{T1}}
\and
\author[B]{\fnms{Carl} \snm{Mueller}\ead[label=e2]{carl.2014@outlook.com}\thanksref{T2}}
\thankstext{T1}{Supported in part by a grant from the Swiss National Foundation for Scientific Research.}
\thankstext{T2}{Supported in part by an NSF grant.}
\address[A]{Institut de Math\'ematiques\\
Ecole Polytechnique F\'ed\'erale de Lausanne\\
Station 8\\
CH-1015 Lausanne\\
Switzerland\\
\printead{e1}}
\affiliation{Ecole Polytechnique F\'ed\'erale de Lausanne and
University of Rochester}
\address[B]{Department of Mathematics\\
University of Rochester\\
Rochester, New York 14627\\
USA\\
\printead{e2}}

\runauthor{R.~C. Dalang and C. Mueller}
\end{aug}

\received{\smonth{2} \syear{2013}}
\revised{\smonth{11} \syear{2013}}

%
\begin{abstract}
It is well known that an $N$-parameter $d$-dimensional Brownian sheet
has no $k$-multiple points when $(k-1)d > 2kN$, and does have such
points when $(k-1)d < 2kN$. We complete the study of the existence of
$k$-multiple points by showing that in the critical cases where $(k-1)
d = 2kN$, there are a.s. no $k$-multiple points.
\end{abstract}

%
\begin{keyword}[class=MSC]
\kwd[Primary ]{60G17}
\kwd[; secondary ]{60G15}
\kwd{60G60}
\end{keyword}

\begin{keyword}
\kwd{Brownian sheet}
\kwd{multiple points}
\kwd{Girsanov's theorem}
\end{keyword}
%
\end{frontmatter}

\section{Introduction and main theorems}
\label{section:statement-of-theorems}
Let $d$ and $N$ be positive integers, and let $B=(B^1,\ldots,B^d)$
denote an $N$-parameter Brownian sheet with values in $\IR^d$, that is,
$B$ is a centered $\IR^d$-valued Gaussian random field with continuous
sample paths, defined on a probability space $(\Omega, \ImF,P)$, with
parameter set $\IR_+^N$ and covariances
\[
\operatorname{Cov}\bigl(B^i(\mathbf{s}), B^j(\mathbf{t})
\bigr) = \delta_{i,j} \prod_{\ell= 1}^N
(s_\ell\wedge t_\ell),
\]
where $\delta_{i,j} = 1$ if $i=j$ and $\delta_{i,j} = 0$ otherwise,
$\mathbf{s}, \mathbf{t} \in\IR_+^N$, $\mathbf{s} = (s_1,\ldots,s_N)$
and $\mathbf{t} = (t_1,\ldots,t_N)$.

The Brownian sheet is perhaps the most studied extension to
multiparameter Gaussian processes of classical Brownian motion, to
which it reduces when $N=1$. Khoshnevisan devotes a chapter to this
process in his book \cite{khosh}. The CIME Summer School lectures \cite
{dal1} contain a presentation of the history of the study of this
random field, and its connections to statistics, Markov properties,
level sets, stochastic partial differential equations, potential theory
and Malliavin calculus.

Here, we are interested in a fundamental sample path property of this
random field, namely multiple points, or self-intersections. For
$\omega\in\Omega$ and integers $k \geq2$, a~point $x \in\IR^d$ is a
\emph{$k$-multiple point of $\mbt\mapsto B(\mbt,\omega)$} if there
exist distinct parameters $\mathbf{t}^1,\ldots,\mathbf{t}^k\in\,
]0,\infty[^N$ such that $B(\mathbf{t}^1,\omega) = \cdots= B(\mathbf
{t}^k,\omega) =x$. We denote the (random, possibly empty) set of all
$k$-multiple points of $\mbt\mapsto B(\mbt,\omega)$ by $M_k(\omega)$.
Note that $M_{k+1}(\omega) \subset M_k(\omega)$.

Typically, for $d$ small and $N$ large, the set of $k$-multiple points
is a.s. nonempty, while for $d$ large and $N$ small, $M_k$ is empty
a.s. See \cite{Dal} for the history of this problem in the case of
Brownian motion ($N=1$).

When $N>1$ and $k\geq2$, it was shown in \cite{K97} that $k$-multiple
points exist if $(k-1)d < 2kN$ and do not exist if $(k-1)d > 2kN$. The
critical case $k=2$ and $d=4N$ was handled in \cite{Dal}, where it was
shown, via quantitative estimates on the conditional distribution of a
pinned Brownian sheet and a decoupling method, that there are no double
points in the critical case. It is also shown in \cite{Dal} that there
are no $k$-multiple points that arise from ordered configurations of
distinct parameters, such as $\mathbf{t}^1 \prec\cdots\prec\mathbf
{t}^k$, where ``$\prec$'' denotes the componentwise (partial) order.

In this paper, we solve the remaining critical cases, where $N > 1$,
$k\geq2$ and $(k-1)d = 2kN$, without any constraints on the parameters
$\mathbf{t}^1,\ldots, \mathbf{t}^k$. The main result of this paper is
the following statement concerning the absence of $k$-multiple points
in these critical cases.

\begin{theorem}
\label{th:0}
Fix $N>1$ and $k\geq2$. If $N$, $d$ and $k$ are such that $(k-1) d =
2kN$, then an $N$-parameter $d$-dimensional Brownian sheet has no
$k$-multiple points, that is, $P\{M_k \neq\varnothing\} =0$.
\end{theorem}

The proof of this theorem relies on known results for hitting
probabilities of the Brownian sheet, due to Khoshnevisan and Shi \cite
{KhoshShi}, on results for intersections of $k$ independent Brownian
sheets, due to Peres \cite{peres}, and a decoupling idea. While \cite{Dal}
used quantitative estimates to obtain their decoupling, 
we will achieve our decoupling here by using Girsanov's theorem. Our
decoupling result is the following.


Let $\mathcal{T}_N^k$ denote the set of parameters $(\mathbf{t}^1,\ldots
,\mathbf{t}^k)$ with $\mathbf{t}^i\in\, ]0,\infty[^N$ such that no two
$\mathbf{t}^i$ and $\mathbf{t}^j$ ($i\ne j$) share a common coordinate:
\begin{eqnarray*}
&&\mathcal{T}_N^k = \bigl\{\bigl(\mathbf{t}^1,
\ldots,\mathbf{t}^k\bigr) \in \bigl(]0,\infty [^N
\bigr)^k\dvtx t^i_\ell\ne t^j_\ell,
\mbox{ for all } \ell=1,\ldots,N
\\
&&\hspace*{192pt} \mbox{and } 1\le i < j\le k\bigr\}
\end{eqnarray*}
[here, $\mbt^i = (t^i_1,\ldots,t^i_N)$, so in our notation, the
coordinates $t^i_\ell$ of $\mbt^i$ inherit the superscript].


%
\begin{theorem}
\label{th:1}
Let $A\subset\IR^d$ be a Borel set. For all $k\in\{2,3,\ldots\}$, we have
\[
P\bigl\{\exists\bigl(\mathbf{t}^1,\ldots,\mathbf{t}^k
\bigr)\in\mathcal{T}_N^k\dvtx B\bigl(\mathbf{t}^1
\bigr)=\cdots=B\bigl(\mathbf{t}^k\bigr)\in A\bigr\}>0
\]
if and only if
\[
P\bigl\{\exists\bigl(\mathbf{t}^1,\ldots,\mathbf{t}^k
\bigr)\in\mathcal{T}_N^k\dvtx W_1\bigl(
\mathbf{t}^1\bigr)=\cdots=W_k\bigl(\mathbf{t}^k
\bigr)\in A\bigr\}>0,
\]
where $W_1,\ldots,W_k$ are independent $N$-parameter Brownian sheets
with values in~$\IR^d$.
\end{theorem}

The proof of this theorem uses an explicit formula for the conditional
expectation $\tilde B(t)$ of $B(t)$ given the values of the sheet in a
product of $N-1$ complements of intervals and a single interval (see
Lemma~\ref{lem3.1}), together with the fact that Girsanov's theorem can
be used to show that the law of the process $B(t) - \tilde B(t)$ is
mutually absolutely continuous with respect to the law of $B$ (see
Lemma~\ref{lem3.5}).

In order to deal with the possibility of a $k$-multiple point arising
from parameters $\mathbf{t}^1,\ldots,\mathbf{t}^k$ that share a common
coordinate, define
\[
\mathcal{H}^k_N(i,j;\ell)= \bigl\{\bigl(
\mathbf{t}^1,\ldots,\mathbf{t}^k\bigr)\in \bigl(]0,
\infty[^N\bigr)^k\dvtx t^i_\ell=t^j_\ell
\bigr\}.
\]
That is, $\mathcal{H}^k_N(i,j;\ell)$ is the set of $(\mathbf{t}^1,\ldots
,\mathbf{t}^k)$ for which $\mathbf{t}^i$ and $\mathbf{t}^j$ share their
$\ell$th coordinate.

Our next theorem states that in the critical case $(k-1)d=2kN$, there
are (with
probability one) no $k$-multiple points arising from parameters in
$\mathcal{H}^k_N(i,j;\ell)$.

\begin{theorem}\label{th:multiple-points-shared-indices}
Suppose $(k-1)d=2kN$, $1\leq i<j\leq k$
and $1\leq\ell\leq N$. Then
\[
P\bigl\{\exists\bigl(\mathbf{t}^1,\ldots,\mathbf{t}^k
\bigr)\in\mathcal{H}^k_N(i,j;\ell)\dvtx B\bigl(
\mathbf{t}^1\bigr)=\cdots=B\bigl(\mathbf{t}^k\bigr)\bigr
\}=0.
\]
\end{theorem}

This theorem is proved by using a covering argument. It requires
checking that certain finite-dimensional distributions of increments of
the Brownian sheet have a uniformly bounded density, provided the
increments are taken at points that are at least $\delta$ units apart
($\delta>0$); see Lemma~\ref{lem2.4}. This uses an explicit formula
for the conditional expectation $\bar B(t)$ of $B(t)$ given the values
of the sheet in a product of $N$ complements of intervals (see Lemma~\ref{lem2.1}).

The paper is structured as follows. First, in Section~\ref{section:proof1}, assuming Theorems \ref{th:1} and \ref
{th:multiple-points-shared-indices}, we easily deduce Theorem~\ref
{th:0} from the results of Khoshnevisan and Shi~\cite{KhoshShi} and
Peres \cite{peres}. Then we prove Theorem~\ref
{th:multiple-points-shared-indices} via an argument based on Hausdorff
dimension, as just mentioned. Finally, in Section~\ref{section:proof2},
we show how to use Girsanov's theorem in order to prove Theorem~\ref{th:1}.

\section{Proof of Theorems \texorpdfstring{\protect\ref{th:0}}{1.1} and \texorpdfstring{\protect\ref{th:multiple-points-shared-indices}}{1.3}}
\label{section:proof1}

We first prove Theorem~\ref{th:0}, assuming Theorems \ref{th:1} and
\ref{th:multiple-points-shared-indices}.

\begin{pf*}{Proof of Theorem~\ref{th:0}} Clearly,
\begin{eqnarray*}
&& P\{ M_k \neq\varnothing\}
\\
& &\qquad\leq P\bigl\{\exists\bigl(\mathbf{t}^1,\ldots,
\mathbf{t}^k\bigr)\in\mathcal{T}_N^k\dvtx B
\bigl(\mathbf{t}^1\bigr)=\cdots=B\bigl(\mathbf{t}^k\bigr)
\bigr\}
\\
&&\qquad\quad{} + \sum_{i=1}^{k-1} \sum
_{j=i+1}^{k}\sum_{\ell=1}^{N}
P\bigl\{ \exists\bigl(\mathbf{t}^1,\ldots,\mathbf{t}^k
\bigr)\in\mathcal{H}^k_N(i,j;\ell)\dvtx 
B\bigl(
\mathbf{t}^1\bigr)=\cdots=B\bigl(\mathbf{t}^k\bigr)\bigr
\}.
\end{eqnarray*}
By Theorem~\ref{th:multiple-points-shared-indices}, the second term
vanishes, and by Theorem~\ref{th:1}, the first term vanishes if and
only if
%
\begin{equation}
\label{intBS} P\bigl\{\exists\bigl(\mathbf{t}^1,\ldots,
\mathbf{t}^k\bigr)\in\mathcal {T}_N^k\dvtx
W_1\bigl(\mathbf{t}^1\bigr)=\cdots=W_k\bigl(
\mathbf{t}^k\bigr)\bigr\} =0,
\end{equation}
where $W_1,\ldots,W_k$ are independent $N$-parameter Brownian sheets
with values in~$\mathbf{R}^d$. According to \cite{KhoshShi}, for all sets of the form
$R = \prod_{\ell=1}^N [s^0_\ell,s^1_\ell] \subset \,]0,\infty[^N$,
there is a finite constant $C \geq1$ such that for all nonrandom Borel
sets $A \subset\IR^d$ contained in a fixed compact subset of $\IR^d$,
\[
C^{-1} \operatorname{Cap}_{d-2N}(A) \leq P\bigl\{\exists\mbt
\in R\dvtx W^i(\mbt) \in A \bigr\} \leq C \operatorname{Cap}_{d-2N}(A),
\]
where $\operatorname{Cap}(\cdot)$ denotes Bessel--Riesz capacity. We recall that
$\operatorname{Cap}(A)$ is defined as follows. Let $\mathcal{P}(K)$ denote the
collection of all probability measures that are supported by the Borel
set $K\subseteq\IR^d$, and define the $\beta$-dimensional \emph
{capacity} of $A$ by
\[
\operatorname{Cap}_\beta(A):= \Bigl[\mathop{\inf_{\mu\in\mathcal{P}(K)\dvtx }}_{
K\subset A\ \mathrm{is\ compact}}
\mathrm{I}_\beta(\mu) \Bigr]^{-1},
\]
where $\inf\varnothing:=\infty$, and $\mathrm{I}_\beta(\mu)$ is
the $\beta$-dimensional \emph{energy} of $\mu$, defined
as follows for all $\mu\in\mathcal{P}(\IR^d)$ and $\beta\in\IR$:
\[
\mathrm{I}_\beta(\mu):= \iint\kappa_\beta(x-y) \mu(d x)
\mu(d y).
\]
In this formula, the function $\kappa_\beta\dvtx \IR^d\to\IR_+\cup\{\infty\}
$ is defined by
\[
\kappa_\beta(x):= %
\cases{ \|x
\|^{-\beta},&\quad $\mbox{if $\beta>0$}$,\vspace*{2pt}
\cr
\log_+\bigl(\|x
\|^{-1}\bigr),&\quad $\mbox{if $\beta=0$}$,\vspace*{2pt}
\cr
1,&\quad $\mbox{if $\beta<0$}$,}
\]
where, as usual, $1/0:=\infty$ and $\log_+(z):=1\vee\log(z)$
for all $z\ge0$.

Since $d-2N>0$ because $(k-1)d = 2kN$, it follows from \cite{peres},
Corollary~15.4, that \eqref{intBS} will hold provided $\operatorname{Cap}_{k(d-2N)}(\IR
^d) =0$. According to \cite{khosh}, Appendix C, Corollary~2.3.1, this
is indeed the case since $k(d-2N) = d$, because we are in the critical
dimension where $(k-1)d = 2kN$.
\end{pf*}


Before proving Theorem~\ref{th:multiple-points-shared-indices}, we
need some preliminary lemmas. For $U \subset\IR_{+}^{N}$, we set
$\mathcal{F}(U) = \sigma(B(\mathbf{t}),  \mathbf{t}\in U)$.

\begin{lemma}\label{lem2.1}
For $\ell= 1, \ldots, N$, fix $0< s_{\ell}^{0} < s^{1}_{\ell}$, and set
\[
R= \prod^{N}_{\ell=1} \bigl[s^{0}_{\ell}
, s^{1}_{\ell}\bigr] \quad\mbox{and} \quad S= \prod
^{N}_{\ell=1} \bigl]s^{0}_{\ell},
s^{1}_{\ell}\bigr[^{c}.
\]
Let $\ImJ$ denote the set of functions from $\{1, \ldots, N\}$ into $\{
0,1\}$. Then for $\mathbf{t} \in R$, set
%
\begin{equation}\qquad
\label{eq:1} \bar{B} (\mathbf{t}) = \sum_{\gamma\in\ImJ}
\biggl(\prod_{\ell\in
\gamma^{-1}(\{1\})} \frac{t_{\ell}- s^{0}_{\ell}}{s^{1}_{\ell}-
s^{0}_{\ell}} \biggr) \biggl(
\prod_{\ell\in\gamma^{-1}(\{0\})} \frac
{s^{1}_{\ell}- t_{\ell}}{s^{1}_{\ell}- s^{0}_{\ell}} \biggr)B
\bigl(s^{\gamma
(1)}_{1}, \ldots, s_{N}^{\gamma(N)}
\bigr)
\end{equation}
(we use the convention that a product over an empty set of indices is
equal to $1$). Then $\bar{B} (\mathbf{t})= E(B(\mathbf{t}) \mid\mathcal
{F}(S))$.
\end{lemma}

\begin{remark} The set of corners (extreme points) of $R$ is
\[
\mathcal{C} = \bigl\{\bigl(s_{1}^{\gamma(1)}, \ldots,
s_{N}^{\gamma(N)}\bigr)\dvtx \gamma \in\ImJ\bigr\},
\]
so the sum over $\gamma$ in \eqref{eq:1} involves $B$ evaluated at each
corner of $R$.
\end{remark}

\begin{pf*}{Proof of Lemma~\ref{lem2.1}} Since the components of $B$
are independent, we may and will assume in this proof that $d=1$. In
this case, since we are working with Gaussian random variables, it
suffices to prove that for each $\mathbf{s}\in S$,
%
\begin{equation}
\label{eq:2} E\bigl(\bar{B}(\mathbf{t}) B(\mathbf{s})\bigr) = E\bigl(B(
\mathbf{t}) B(\mathbf{s})\bigr).
\end{equation}
The right-hand side of \eqref{eq:2} is equal to $\prod^{N}_{\ell=1}
(t_{\ell} \wedge s_{\ell})$, so we compute the left-hand side of \eqref
{eq:2}. Clearly,
\begin{eqnarray*}
E\bigl(\bar{B} (\mathbf{t})B(\mathbf{s})\bigr) &=& \sum
_{\gamma\in\ImJ} \biggl(\prod_{\ell\in\gamma^{-1}(\{1\})}
\frac{t_{\ell}- s^{0}_{\ell
}}{s^{1}_{\ell}- s^{0}_{\ell}} \biggr) \biggl(\prod_{\ell\in\gamma^{-1}(\{0\})}
\frac{s^{1}_{\ell}- t_{\ell
}}{s^{1}_{\ell}- s^{0}_{\ell}} \biggr) \prod^{N}_{\ell=1}
\bigl(s_{\ell}^{\gamma(\ell)}\wedge s_{\ell}\bigr)
\\
& =& \prod^{N}_{\ell=1} \biggl[
\bigl(s^{1}_{\ell} \wedge s_{\ell}\bigr)
\frac
{t_{\ell} - s^{0}_{\ell}}{s^{1}_{\ell}-s^{0}_{\ell}} + \bigl(s^{0}_{\ell} \wedge s_{\ell}
\bigr) \frac{s^{1}_{\ell}-t_{\ell}}{s^{1}_{\ell}-s^{0}_{\ell
}} \biggr].
\end{eqnarray*}
Therefore, \eqref{eq:2} will be proved if we show that for each $\ell
\in\{1, \ldots, N\}$,
%
\begin{equation}
\label{eq:3} t_{\ell} \wedge s_{\ell} = \bigl(s^{1}_{\ell}
\wedge s_{\ell}\bigr) \frac{t_{\ell
}-s^{0}_{\ell}}{s^{1}_{\ell}-s^{0}_{\ell}} + \bigl(s^{0}_{\ell}
\wedge s_{\ell}\bigr) \frac{s^{1}_{\ell}-t_{\ell}}{s^{1}_{\ell}-s^{0}_{\ell}}.
\end{equation}
There are two cases to distinguish.

\emph{Case} 1. $s_{\ell} \leq s^{0}_{\ell}$. In this case,
$s^{k}_{\ell} \wedge s_{\ell} = s_{\ell}$ for $k \in\{0,1\}$ and
$t_{\ell} \wedge s_{\ell} = s_{\ell}$, since $s^{0}_{\ell}\leq t_{\ell}
\leq s_{\ell}^{1}$, so the right-hand side of \eqref{eq:3} is equal to
\[
s_{\ell} \frac{t_{\ell} - s^{0}_{\ell}}{s^{1}_{\ell}-s^{0}_{\ell}}+ s_{\ell} \frac{s^{1}_{\ell} - t_{\ell}}{s^{1}_{\ell}-s^{0}_{\ell}}=
s_{\ell},
\]
which is also the left-hand side of \eqref{eq:3}.

\emph{Case} 2. $s_{\ell} \geq s^{1}_{\ell}$. In this case,
$s^{k}_{\ell}\wedge s_{\ell} = s^{k}_{\ell}$ for $k \in\{0,1\}$ and
$t_{\ell} \wedge s_{\ell} = t_{\ell}$, so the right-hand side of \eqref
{eq:3} is equal to
\[
s^{1}_{\ell} \frac{t_{\ell} - s^{0}_{\ell}}{s^{1}_{\ell}-s^{0}_{\ell
}}+ s^{0}_{\ell}
\frac{s^{1}_{\ell} - t_{\ell}}{s^{1}_{\ell}-s^{0}_{\ell
}}= t_{\ell},
\]
and which is also the left-hand side of \eqref{eq:3}.

This completes the proof of Lemma~\ref{lem2.1}.
\end{pf*}

\begin{remark}\label{rem2.3} We note that the right-hand side of \eqref{eq:1} is in
fact a convex combination of the values of $B$ at the corners of $R$,
since each coefficient is nonnegative and
\begin{eqnarray*}
&&\sum_{\gamma\in\mathcal{J}} \biggl( \prod
_{\ell\in\gamma^{-1} (\{1\}
)}\frac{t_{\ell} - s^{0}_{\ell}}{s^{1}_{\ell}-s^{0}_{\ell}} \biggr) \biggl( \prod
_{\ell\in\gamma^{-1} (\{0\})}\frac{s^{1}_{\ell} - t_{\ell
}}{s^{1}_{\ell}-s^{0}_{\ell}} \biggr)
\\
&&\qquad = \prod^{N}_{\ell=1} \biggl[
\frac{t_{\ell} - s^{0}_{\ell
}}{s^{1}_{\ell}-s^{0}_{\ell}} + \frac{s^{1}_{\ell} - t_{\ell
}}{s^{1}_{\ell}-s^{0}_{\ell}} \biggr]=1.
\end{eqnarray*}
\end{remark}

\begin{lemma}\label{lem2.4} Fix $\delta>0$ (small), $K \in\IN$ (positive and
large), and $k \in\IN$, $k \geq2$.
\begin{longlist}[(a)]
\item[(a)] There is $C >0$ such that for all $\mathbf{t}^{1}, \ldots, \mathbf
{t}^{k}$ such that $\Vert\mathbf{t}^{i} - \mathbf{t}^{j}\Vert\geq
\delta$, for all $i \neq j$ with $i, j \in\{1, \ldots, k\}$, and $K
\geq t^{i}_{\ell} \geq\delta$, for all $\ell= 1, \ldots, N$ and $i \in
\{1, \ldots, k\}$, the random vector $(B(\mathbf{t}^{1}), \ldots, B(\mathbf
{t}^{k}))$ has a joint probability density function that is bounded by $C$.

\item[(b)] For the same choices of $\mathbf{t}^{1}, \ldots, \mathbf{t}^{k}$,
the ${(\IR^{d})}^{k-1}$-valued random vector
\[
\bigl(B\bigl(\mathbf{t}^{1}\bigr)- B\bigl(\mathbf{t}^{2}
\bigr),B\bigl(\mathbf{t}^{2}\bigr)-B\bigl(\mathbf {t}^{3}
\bigr), \ldots, B\bigl(\mathbf{t}^{k-1}\bigr)-B\bigl(
\mathbf{t}^{k}\bigr)\bigr)
\]
has a bounded probability density function (with bound depending only
on $\delta$, $K$ and $k$, as well as $d$ and $N$).
\end{longlist}
\end{lemma}

\begin{pf}Since the $B^{1}, \ldots, B^{d}$ are independent Brownian sheets,
we may and will assume in this proof that $d=1$.

We first deduce (b) from (a). Let
\[
Y= \bigl(B\bigl(\mathbf{t}^{1}\bigr)- B\bigl(\mathbf{t}^{2}
\bigr), \ldots, B\bigl(\mathbf{t}^{k-1}\bigr) -B\bigl(
\mathbf{t}^{k}\bigr), B\bigl(\mathbf{t}^{k}\bigr)\bigr).
\]
Then $Y$ is obtained from $(B(\mathbf{t}^{1}), \ldots, B(\mathbf
{t}^{k}))$ by applying an invertible linear transformation from ${(\IR
^{d})}^{k}$ into ${(\IR^{d})}^{k}$. Therefore, by (a), $Y$ has a
bounded joint probability density function. It follows that the
probability density function of $(B(\mathbf{t}^{1})- B(\mathbf{t}^{2}),
\ldots, B(\mathbf{t}^{k-1}) -B(\mathbf{t}^{k}))$, which is a marginal
density of $Y$, is bounded by the same constant. This proves (b).

We now prove (a). Set
\[
n = \inf \biggl\{n \in\IN\dvtx 2^{-n} < \frac{\delta}{3 \sqrt{N}} \biggr\},
\]
and consider a dyadic grid in $\IR_{+}^{N}$ with edges of length
$2^{-n}$. We let $G_{\delta,K}$ denote the set of such grid points with
all coordinates $\leq K$.

By construction, each closed box in this grid contains at most one of
the $\mathbf{t}^{i}$, and we denote by $R^{i}$ the box containing
$\mathbf{t}^{i} $. Suppose that
\[
R^{i} = \prod^{N}_{\ell=1}
\bigl[s_{\ell}^{i, 0}, s_{\ell}^{i,1}\bigr]\quad
\mbox{and set}\quad S^{i} = \prod^{N}_{\ell=1}
\bigl]s_{\ell}^{i, 0}, s_{\ell}^{i,1}\bigr[^c.
\]
Because of our choice of $n$, the set $\ImC^{i}$ of corners of $R^{i}$
is distinct from $\ImC^{j}$ when $i\neq j$.

Define
\[
Y^{i} = E\bigl(B \bigl(\mathbf{t}^{i}\bigr) \vert
\mathcal{F} \bigl(S^{i}\bigr)\bigr), \qquad i= 1, \ldots, k.
\]
Then $B(\mathbf{t}^{i})- Y^{i}$ is orthogonal to $Y^{i}$, and for $j\neq i$, since $Y^{j}$ is a linear combination of values of $B$ at
elements of $S^{i}$ (because $\ImC^{j} \cap\ImC^{i}= \varnothing$),
$B(\mathbf{t}^{i})- Y^{i}$ is orthogonal to $Y^{j}$. Letting $Y =
(Y^{1}, \ldots, Y^{k})$ and $Z = (B(\mathbf{t}^{1})- Y^{1}, \ldots,
B(\mathbf{t}^{k})- Y^{k})$, we see that the Gaussian vectors $Y$ and
$Z$ are independent, and
\[
\bigl(B\bigl(\mathbf{t}^{1}\bigr), \ldots, B\bigl(
\mathbf{t}^{k}\bigr)\bigr)= Y + Z.
\]

Using properties of convolution, we see that it suffices to
show that the joint probability density function of $Y$ is bounded
[uniformly over the $(\mathbf{t}^{1}, \ldots, \mathbf{t}^{k})$].

Since $Y$ is a Gaussian random vector, let $M$ be its
variance--covariance matrix. It suffices to show that
%
\begin{equation}
\label{eq:4} \det M > c >0,
\end{equation}
where $c$ depends only on $\delta$, $K$ and $k$, as well as $d$ and $N$.

Consider the random vector $(B(\mathbf{r}),  \mathbf{r} \in G_{\delta
,K})$. Observe that this random vector can be obtained by applying an
invertible linear transformation, from $\IR^{({(2^{n}K)}^{N})}$ into
itself (recall that $d=1$), to the random vector ($W(R)$, $R$ a box in
the grid), which has i.i.d. components, each with variance
${(2^{-n})}^{N} >0$. Therefore, $(B(\mathbf{r}), \mathbf{r }\in
G_{\delta,K})$ has a bounded density, where the bound depends only on
$\delta$ and $K$ (and $d$ and $N$). This implies that $(B( \mathbf
{t}),  t \in\ImC^{i},  i= 1, \ldots, k)$ has a joint probability
density function that is bounded, since it is a marginal density of
$(B( \mathbf{r}),  \mathbf{r} \in G_{\delta,K})$.

Let $\tilde{M}$ be the variance--covariance matrix of the Gaussian
random vector $(B(\mathbf{t}),  \mbt\in\ImC^{i},  i= 1, \ldots, k)$.
Then by the above, there is $c>0$ such that $\det\tilde{M} > C$. In
particular, there is $c_{0} > 0$ such that
\[
\lambda^{T} \tilde{M} \lambda\geq c_{0} \Vert\lambda
\Vert^{2}\qquad \mbox{for all } \lambda\in\IR^{k2^{N}}.
\]
Note that $c_{0}$ depends only on $(\delta, K, k, d, N)$.

Let $\mu\in\IR^{k}$. Then
\begin{eqnarray*}
\mu^{T} M \mu & = & \operatorname{Var} \Biggl( \sum
_{i=1}^{k} \mu_{i} Y^{i}
\Biggr)
\\
& = & \operatorname{Var} \Biggl( \sum_{i=1}^{k}
\mu_{i} \sum_{\mbs^{i,j}\in
\ImC^{i}} a_{i,j} B\bigl(
\mathbf{s}^{i,j}\bigr) \Biggr)
\\
& \geq& c_{0} \sum_{i=1}^{k}
\sum_{\mbs^{i,j}\in\ImC^{i}}\mu_{i}^{2}
a_{i,j}^{2},
\end{eqnarray*}
where the $a_{i,j}$ are the coefficients obtained in formula \eqref
{eq:1} of Lemma~\ref{lem2.1}. According to Remark~\ref{rem2.3}, $\sum_{s^{i,j}\in\ImC^{i}} a_{i,j}=1$ and $a_{i,j} \geq0$, therefore,
there is $\alpha>0$ such that $ \sum_{s^{i,j}\in\ImC^{i}} a_{i,j}^{2}
> \alpha$. We conclude that
\[
\mu^{T} M \mu\geq c_{0} \alpha\sum
^{k}_{i=1} \mu_{i}^{2},
\]
and this implies that det $M > c_{1} >0$, where $c_{1}$ depends only on
$(\delta, K, k, d, N)$. In turn, this proves \eqref{eq:4} and completes
the proof of (a) in Lemma~\ref{lem2.4}.
\end{pf}

\begin{pf*}{Proof of Theorem~\ref
{th:multiple-points-shared-indices}} It suffices to prove the theorem
in the case where $i=1, j=2$ and $\ell=1$. Therefore, we write
$\mathcal{H}_{N}^{k}$ instead of $\mathcal{H}_{N}^{k}(1, 2; 1)$.

For $\delta>0$, set
\begin{eqnarray*}
&&\mathcal{H}_{N}^{k}(\delta) = \bigl\{\bigl(
\mathbf{t}^{1}, \ldots, \mathbf {t}^{k}\bigr) \in
\mathcal{H}_{N}^{k}\dvtx t_{\ell}^{i}
\geq\delta, \bigl\Vert \mathbf{t}^{i}- \mathbf{t}^{j}\bigr\Vert\geq
\delta,
\\
&&\hspace*{50pt} \mbox{for all } i \neq j, \ell=1, \ldots, N, i, j \in\{1, \ldots, k\} \bigr
\}.
\end{eqnarray*}
Since $\mathcal{H}_{N}^{k} = \bigcup_{n=1}^{\infty} \mathcal{H}_{N}^{k}
( \frac{1}{n}  )$, it suffices to prove that for fixed
$\delta> 0$,
\[
P\bigl\{\exists\bigl(\mathbf{t}^{1}, \ldots, \mathbf{t}^{k}
\bigr) \in\mathcal {H}_{N}^{k}(\delta) \dvtx B \bigl(
\mathbf{t}^{1}\bigr)= \cdots= B\bigl(\mathbf{t}^{k}\bigr)
\bigr\} =0.
\]

Consider the random field indexed by $(]0,\infty[^N)^k$ with values in
${(\IR^{d})}^{k-1}$ defined by
\[
X \bigl(\mathbf{t}^{1}, \ldots, \mathbf{t}^{k}\bigr)=
\bigl(B \bigl(\mathbf{t}^{1}\bigr)- B\bigl(\mathbf{t}^{2}
\bigr), B \bigl(\mathbf{t}^{2}\bigr)- B\bigl(\mathbf{t}^{3}
\bigr), \ldots, B \bigl(\mathbf{t}^{k-1}\bigr)- B\bigl(
\mathbf{t}^{k}\bigr)\bigr).
\]
Then
\[
B \bigl(\mathbf{t}^{1}\bigr)= \cdots= B\bigl(\mathbf{t}^{k}
\bigr)\quad\iff\quad X\bigl( \mathbf{t}^{1}, \ldots, \mathbf{t}^{k}
\bigr) =0,
\]
so parameters which give rise to a $k$-multiple point of $B$ are
$k$-tuples at which $X$ hits $0$ ($\in{(\IR^{d})}^{k-1}$). Therefore,
it will suffice to show that
%
\begin{equation}
\label{eq:5} P\bigl\{\exists\bigl(\mathbf{t}^{1}, \ldots,
\mathbf{t}^{k}\bigr) \in\mathcal {H}_{N}^{k} (
\delta)\dvtx X\bigl(\mathbf{t}^{1}, \ldots, \mathbf{t}^{k}
\bigr)=0\bigr\} =0.
\end{equation}

Let $D(K) = \mathcal{H}_{N}^{k} (\delta) \cap{([0,K]^{N})}^{k}$.
Since $\mathcal{H}_{N}^{k}$ is a vector space of dimension $kN-1$,
there is $C>0$ such that for all large $n \geq1$, we can cover $D(K)$
by $C{(2^{2n})}^{kN-1}$ dyadic boxes in $(\IR^N)^k$ with edges of
length $2^{-2n}$. Let $\ImD_{n}$ be the set of boxes in such a
covering, and for $D \in\ImD_{n}$, let $t_{n}(D)$ be the corner of $D$
for which all coordinates are smallest possible.

For $(\mathbf{t}^{1}, \ldots, \mathbf{t}^{k}) \in D$, let $p_{(\mathbf
{t}^{1}, \ldots, \mathbf{t}^{k})} (z_{1}, \ldots, z_{k-1})$ be the value
of the joint probability density function of $X(\mathbf{t}^{1}, \ldots,
\mathbf{t}^{k})$ at $(z_{1}, \ldots, z_{k-1})\in{(\IR^{d})}^{k-1}$. By
Lemma~\ref{lem2.4}, there is $C< +\infty$ such that
%
\begin{equation}
\label{eq:6} p_{(\mathbf{t}^{1}, \ldots, \mathbf{t}^{k})} (z_{1}, \ldots, z_{k-1})\leq
C.
\end{equation}
Let $B(0, n2^{-n})$ denote the ball in ${(\IR^{d})}^{k-1}$ centered at
0 with radius $n 2^{-n}$. By \eqref{eq:6},
%
\begin{equation}
\label{eq:7} P\bigl\{X\bigl(\mathbf{t}^{1}, \ldots,
\mathbf{t}^{k}\bigr)\in B\bigl(0, n2^{-n}\bigr)\bigr\} \leq C
{\bigl(n 2^{-n}\bigr)}^{d(k-1)}.
\end{equation}
In order to prove \eqref{eq:5}, it suffices to prove \eqref{eq:5} with
$\mathcal{H}_{N}^{k}(\delta)$ replaced by $D(K)$. So, we compute
\begin{eqnarray*}
&&P\bigl\{ \exists\bigl(\mathbf{t}^{1}, \ldots, \mathbf{t}^{k}
\bigr) \in D(K)\dvtx X\bigl(\mathbf {t}^{1}, \ldots,
\mathbf{t}^{k}\bigr)=0\bigr\}
\\
&&\qquad \leq P\bigl\{\exists\bigl(\mathbf{t}^{1}, \ldots,
\mathbf{t}^{k}\bigr) \in D(K)\dvtx X \bigl(\mathbf{t}^{1},
\ldots, \mathbf{t}^{k}\bigr) \in B\bigl(0, 2^{-n}\bigr)\bigr
\}
\\
&&\qquad \leq \sum_{D \in\ImD_{n}} P\bigl\{\exists\bigl(
\mathbf{t}^{1}, \ldots, \mathbf{t}^{k}\bigr) \in D\dvtx X
\bigl(\mathbf{t}^{1}, \ldots, \mathbf{t}^{k}\bigr) \in B
\bigl(0, 2^{-n}\bigr)\bigr\}
\\
&&\qquad \leq \sum_{D \in\ImD_{n}} P \Bigl(\bigl\{X
\bigl(t_{n}(D)\bigr) \in B\bigl(0, n 2^{-n}\bigr)\bigr\}
\\
&&\hspace*{39pt}{}\qquad\quad \cup \Bigl\{ \sup_{t\in D}\bigl\Vert X(t) -X\bigl(t_{n}(D)
\bigr) \bigr\Vert\geq(n-1) 2^{n} \Bigr\} \Bigr).
\end{eqnarray*}
We now use \eqref{eq:7} to bound this by
\begin{eqnarray*}
2^{2n(kN-1)}  \Bigl[C{\bigl(n2^{-n}\bigr)}^{d(k-1)}
 + \sup_{D \in\ImD_{n}} P \Bigl\{\sup_{t\in D}\bigl\Vert X(t)
-X\bigl(t_{n}(D)\bigr) \bigr\Vert\geq(n-1) 2^{-n} \Bigr\} \Bigr].
\end{eqnarray*}
It follows from the scaling property of the Brownian sheet (\cite{walsh},
Chapter~1) that the supremum over $D \in\ImD_{n}$ is no
greater than that achieved by the box $D^* = [K-2^{-2n},K]^{Nk}$, and
we will show below that
%
\begin{equation}
\label{eq:8}
\lim_{n\to+\infty} 2^{2n(kN-1)}P \Bigl\{
\sup_{t\in D^*}\bigl\Vert X(t) -X\bigl(t_{n}\bigl(D^*\bigr)
\bigr) \bigr\Vert\geq(n-1) 2^{-n} \Bigr\} = 0, 
\end{equation}
so it remains to examine the term $n^{d(k-1)} {(2^{-n})}^{d(k-1)
-2kN+2}$. Since we are in the critical case, $2kN=(k-1)d$, so the
exponent of $2^{-n}$ is equal to 2 and, therefore,
\[
n^{d(k-1)} {\bigl(2^{-n}\bigr)}^{d(k-1) -2kN+2} =
n^{d(k-1)}2^{-2n} \to0
\]
as $n\to+\infty$. This will prove \eqref{eq:7} and complete the proof
of Theorem~\ref{th:multiple-points-shared-indices} once we establish
\eqref{eq:8}, to which we now turn.

We can write $D^* = D_1 \times\cdots\times D_k$, where each $D_i$ is
a box in $\IR^N$ with edges of length $2^{-2n}$, and we can write
$t_n(D^*) = (t^1_n(D_1),\ldots,t^k_n(D_k))$. Clearly,
\[
\bigl\Vert X(t) -X\bigl(t_{n}\bigl(D^*\bigr)\bigr)\bigr \Vert\leq2 \sum
_{i=1}^k\bigl \Vert B\bigl(\mbt^i
\bigr) - B\bigl(t^i_n(D_i)\bigr) \bigr\Vert,
\]
so it suffices to prove that for each $i\in\{1,\ldots,k\}$ and $n$
sufficiently large, there are constants $C < \infty$ and $c>0$ such that
%
\begin{equation}
\label{eq2.9a} P \biggl\{\sup_{\mbt^i \in D_i} \bigl\Vert B\bigl(
\mbt^i\bigr) - B\bigl(t^i_n(D_i)
\bigr)\bigr \Vert \geq\frac{(n-1)2^{-n}}{2k} \biggr\} \leq C e^{-c^2 (n-1)^2/2}.
\end{equation}

In order to simplify the notation, we assume that $D_i =
[1,1+2^{-2n}]^N$, so $t^i_n(D_i) = (1,\ldots,1)$, and we write $\mbt^i =
(t^i_1,\ldots,t^i_N)$. We use the decomposition of the Brownian sheet
presented in \cite{kendall}, proof of Theorem~(1.1), to write
\[
B\bigl(\mbt^i\bigr) - B\bigl(t^i_n(D_i)
\bigr) = \sum_{m=1}^N \sum
_{1\leq\ell_1<\cdots
<\ell_m\leq N} W^{(m)}_{\ell_1,\ldots,\ell_m}\bigl(t^i_{\ell_1}-1,
\ldots,t^i_{\ell_m}-1\bigr),
\]
where the $W^{(m)}_{\ell_1,\ldots,\ell_m}$ are $m$-parameter Brownian
sheets and all are mutually independent. There are $2^N-1$ terms in
this decomposition, so, using the scaling property of the Brownian
sheet, we see that
\begin{eqnarray*}
&&P \biggl\{\sup_{\mbt^i \in D_i} \bigl\Vert B\bigl(\mbt^i\bigr) -
B\bigl(t^i_n(D_i)\bigr)\bigr \Vert \geq
\frac{(n-1)2^{-n}}{2k} \biggr\}
\\
&&\qquad \leq\sum_{m=1}^N \sum
_{1\leq\ell_1<\cdots<\ell_m\leq N} P \biggl\{\sup_{\mbt\in[0,1]^m}
W^{(m)}_{\ell_1,\ldots,\ell_m}(\mbt) \geq \frac{(n-1)2^{(m-1)n}}{2k2^N} \biggr\}.
\end{eqnarray*}
Using \cite{OP}, Lemma~1.2, we see that the largest probability in this
sum is obtained when $m=1$, and in this case it is bounded by $4^N P\{
Z\geq c(n-1) \}$, where $Z$ is a standard normal random variable and
$c= 2^{-N-1}/k$. Therefore,
\[
P \biggl\{\sup_{\mbt^i \in D_i} \bigl\Vert B\bigl(\mbt^i\bigr) - B
\bigl(t^i_n(D_i)\bigr)\bigr \Vert \geq
\frac{(n-1)2^{-n}}{2k} \biggr\} \leq N! \,8^N e^{-c^2 (n-1)^2/2},
\]
which proves \eqref{eq2.9a} and completes the proof of Theorem~\ref
{th:multiple-points-shared-indices}.
\end{pf*}

\section{Proof of Theorem \texorpdfstring{\protect\ref{th:1}}{1.2}}
\label{section:proof2}

The main ingredient in the proof of Theorem~\ref{th:1} is the following result.

\begin{theorem}\label{thm3.4} Let $W_{1}, \ldots, W_{k}$ be independent Brownian
sheets. Fix $M>0$ and let $\ImR_{M}$ denote the set of $k$-tuples of
boxes $(R_{1}, \ldots, R_{k})$, where each box $R_{i}$ is contained in
${[M^{-1}, M]}^{N}$ and for each coordinate axis, the projections of
the $R_{i}$ onto this coordinate axis are pairwise disjoint. Then, for
all $(R_{1}, \ldots, R_{k}) \in\ImR_{M}$, the random vectors
\[
\bigl(B\vert_{R_{1}}, \ldots, B \vert_{R_{k}}\bigr)\quad \mbox{and}\quad
\bigl(W_{1} \vert_{R_{1}}, \ldots, W_{k}
\vert_{R_{k}}\bigr)
\]
[with values in $(C(R_{1}, \IR^{d})\times\cdots\times C(R_{k}, \IR
^{d}))$] have mutually absolutely continuous probability distributions.
\end{theorem}

\begin{remark} Using the results of Walsh \cite{walsh1} on propagation
of singularities in the Brownian sheet, it is easy to see that the
conclusion of Theorem~\ref{thm3.4} does not remain valid without the
assumption that the projections of the $R_i$ onto each axis are
pairwise disjoint.
\end{remark}

Before proving Theorem~\ref{thm3.4}, we show that it readily implies
Theorem~\ref{th:1}.

\begin{pf*}{Proof of Theorem~\ref{th:1}} Let $A \subset\IR^d$ be a
Borel set. Fix $M>0$ and set $\mathcal{T}_{N}^{k}(M) = \mathcal
{T}_{N}^{k} \cap{[M^{-1}, M]}^{N}$. Then $\mathcal{T}_{N}^{k} = \bigcup
^{\infty}_{M=1} \mathcal{T}_{N}^{k}(M)$. Therefore,
%
\begin{equation}
\label{eq3.3} P\bigl\{\exists\bigl(\mathbf{t}^{1}, \ldots,
\mathbf{t}^{k}\bigr) \in\mathcal {T}_{N}^{k}
\dvtx B\bigl(\mathbf{t}^{1}\bigr) = \cdots= B\bigl(
\mathbf{t}^{k}\bigr) \in A \bigr\} =0
\end{equation}
is equivalent to
\[
\forall M \in\IN^{*}\qquad P\bigl\{\exists\bigl(\mathbf{t}^{1},
\ldots, \mathbf {t}^{k}\bigr) \in\mathcal{T}_{N}^{k}(M)
\dvtx B\bigl(\mathbf{t}^{1}\bigr) = \cdots= B\bigl(
\mathbf{t}^{k}\bigr) \in A \bigr\} =0,
\]
and this in turn is equivalent to
%
\begin{eqnarray}
\label{eq:11} &&\forall M\in\IN^{*}, \forall(R_{1},
\ldots, R_{k}) \in\ImR_{M}
\nonumber
\\[-8pt]
\\[-8pt]
\nonumber
&&\qquad P\bigl\{ \exists\bigl(\mathbf{t}^{1}, \ldots, \mathbf{t}^{k}
\bigr)\in R_{1} \times\cdots\times R_{k}\dvtx B\bigl(
\mathbf{t}^{1}\bigr)= \cdots=B\bigl(\mathbf{t}^{k}\bigr) \in
A\bigr\} =0.
\end{eqnarray}
Similarly, the property
%
\begin{equation}
\label{eq3.5} P\bigl\{\exists\bigl(\mathbf{t}^{1}, \ldots,
\mathbf{t}^{k}\bigr) \in\mathcal {T}_{N}^{k}
\dvtx W_{1}\bigl(\mathbf{t}^{1}\bigr) = \cdots=
W_{k}\bigl(\mathbf{t}^{k}\bigr) \in A \bigr\} =0
\end{equation}
is equivalent to
%
\begin{eqnarray}
\label{eq:12} &&\forall M\in\IN^{*}, \forall(R_{1}, \ldots,
R_{k}) \in\ImR_{M}\dvtx
\nonumber
\\[-8pt]
\\[-8pt]
\nonumber
&&\qquad P\bigl\{ \exists\bigl(\mathbf{t}^{1}, \ldots, \mathbf{t}^{k}
\bigr)\in R_{1} \times\cdots\times R_{k}\dvtx
W_{1}\bigl(\mathbf{t}^{1}\bigr)= \cdots=W_{k}
\bigl(\mathbf {t}^{k}\bigr) \in A\bigr\} =0.
\end{eqnarray}
According to Theorem~\ref{thm3.4}, properties \eqref{eq:11} and \eqref
{eq:12} are equivalent and, therefore, \eqref{eq3.3} and \eqref{eq3.5}
are also equivalent. This proves Theorem~\ref{th:1}.
\end{pf*}

For Theorem~\ref{thm3.4}, we will need a variant of Lemma~\ref{lem2.1}.

\begin{lemma}\label{lem3.1} For $\ell=1, \ldots, N$, fix $0< s_{\ell}^{0} < s_{\ell
}^{1}$ and set
\[
R= \prod^{N}_{\ell=1} \bigl[s_{\ell}^{0},
s_{\ell}^{1}\bigr] \quad\mbox{and}\quad S = \Biggl( \prod
^{N-1}_{\ell=1} {\bigl]s_{\ell}^{0},
s_{\ell}^{1}\bigr[}^{c} \Biggr) \times\bigl[0,
s^{0}_{N}\bigr].
\]
Let $\mathcal{J}_{N}$ denote the set of functions from $\{1, \ldots,
N-1\}$ into $\{0, 1\}$ and set
\[
\ImC_{N} = \bigl\{ \bigl(s_{1}^{\gamma(1)},
s_{2}^{\gamma(2)}, \ldots, s_{N}^{\gamma(N-1)},
s_{N}^{0} \bigr)\dvtx \gamma\in\mathcal{J}_{N}
\bigr\}.
\]
For $\mbt\in R$, set
\begin{eqnarray*}
\tilde{B} (\mbt) &=& \sum_{\gamma\in\mathcal{J}_{N}} \biggl( \prod
_{\ell\in\gamma^{-1} (\{1\})} \frac{t_{\ell} - s_{\ell}^{0}}{s_{\ell
}^{1}-s_{\ell}^{0}} \biggr) \biggl( \prod
_{\ell\in\gamma^{-1} (\{0\})} \frac{s^{1}_{\ell} - t_{\ell}}{s_{\ell}^{1}-s_{\ell}^{0}} \biggr)
\\
&&\hspace*{23pt}{} \times B \bigl(s_{1}^{\gamma(1)},\ldots, s_{N-1}^{\gamma(N-1)},
s_{N}^{0} \bigr).
\end{eqnarray*}
Then $\tilde{B} (\mbt)= E( B(\mbt) \mid\ImF(S))$.
\end{lemma}

\begin{remark} $\ImC_{N}$ is the set of corners of $R$ with the
smallest of the two possible $N$th coordinates, and $S_{N}$ is in the
``past'' of $R$ if we define the ``past'' using the (partial) order
$\mathbf{s} \leq_{N} \mathbf{t}$ if and only if $s_{N} \leq t_{N}$.
\end{remark}

\begin{pf*}{Proof of Lemma~\ref{lem3.1}} Since the components of $B$
are independent, we may and will assume in this proof that $d=1$. In
this case, as in the proof of Lemma~\ref{lem2.1}, it suffices to prove
that for each $\mathbf{s}\in S$,
%
\begin{equation}
\label{eq:9} E(\tilde{B}_{N} (\mathbf{t}) B(\mathbf{s})= E\bigl( B(
\mathbf{t}) B(\mathbf{s})\bigr).
\end{equation}

The right-hand side of \eqref{eq:9} is equal to $s_{N} \prod_{\ell=1}^{N-1} (t_{\ell} \wedge s_{\ell})$, so we compute the
left-hand side of \eqref{eq:9}. Clearly,
\begin{eqnarray*}
E\bigl(\tilde{B} (\mathbf{t}) B(\mathbf{s})\bigr) &=&s_{N} \sum
_{\gamma\in
{\mathcal{J}}_{N}} \biggl[ \prod_{\ell\in\gamma^{-1} (\{1\})}
\frac
{t_{\ell} - s_{\ell}^{0}}{s_{\ell}^{1}-s_{\ell}^{0}} \biggr] \biggl[ \prod_{\ell\in\gamma^{-1} (\{0\})}
\frac{s^{1}_{\ell} - t_{\ell}}{s_{\ell
}^{1}-s_{\ell}^{0}} \biggr]\bigl(s_{\ell}^{\gamma(\ell)} \wedge
s_{\ell}\bigr)
\\
&=& s_{N} \prod_{\ell=1}^{N-1}
\biggl[\bigl(s_{\ell}^{1} \wedge s_{\ell
}\bigr)
\frac{t_{\ell} - s_{\ell}^{0}}{s_{\ell}^{1}-s_{\ell}^{0}} + \bigl(s_{\ell
}^{0} \wedge s_{\ell}
\bigr) \frac{s^{1}_{\ell} - t_{\ell}}{s_{\ell
}^{1}-s_{\ell}^{0}} \biggr],
\end{eqnarray*}
so \eqref{eq:9} will be proved if we check that for each $\ell\in\{1,
\ldots, N-1\}$,
\[
t_{\ell} \wedge s_{\ell} = \bigl(s_{\ell}^{1}
\wedge s_{\ell}\bigr) \frac{t_{\ell
} - s_{\ell}^{0}}{s_{\ell}^{1}-s_{\ell}^{0}} + \bigl(s_{\ell}^{0}
- s_{\ell
}\bigr) \frac{s^{1}_{\ell} - t_{\ell}}{s_{\ell}^{1}-s_{\ell}^{0}}.
\]
But this is simply equality \eqref{eq:3}, and the proof of Lemma~\ref
{lem3.1} is complete.
\end{pf*}

We will need the following form of Girsanov's theorem for the Brownian
sheet, which is essentially the version given in \cite{NP}, Proposition~1.6. Fix $M >0$. Define the one-parameter filtration $\ImG= (\ImG
_{u},  u \in[0, M])$ by
%
\begin{equation}
\label{fG} \ImG_{u} = \sigma \bigl\{ B(t_{1}, \ldots,
t_{N-1}, v) \dvtx (t_{1}, \ldots, t_{N-1}) \in
\IR_{+}^{N-1}, v \in[0,u] \bigr\}
\end{equation}
(the filtration is completed and made right-continuous). Let $(Z( \mbs
),\break  \mbs\in\IR_{+}^{N-1} \times[0, M])$ be a (jointly measurable)
$\IR^d$-valued random field that is adapted to $\ImG$, that is, for all
$\mbs\in\IR_{+}^{N-1} \times[0, M]$, $Z(\mbs)$ is $\ImG
_{s_{N}}$-measurable. Suppose that
%
\begin{equation}
\label{eq:10} E \biggl(\int_{\IR_{+}^{N-1}\times[0, M]} \bigl\Vert Z(\mathbf{s})
\bigr\Vert^{2} \,d\mathbf{s} \biggr) < +\infty.
\end{equation}
For $u \in[0, M]$, define
\[
L_{u} = \exp \biggl( \int_{\IR_{+}^{N-1}\times[0, u]} Z(\mathbf{s})
\cdot d B(\mathbf{s}) - \frac{1}{2} \int_{\IR_{+}^{N-1}\times[0, u]} \bigl\Vert Z(
\mathbf{s})\bigr\Vert^{2} \,d\mathbf{s} \biggr),
\]
where ``$\cdot$'' denotes the Euclidean inner product and, for each
component, the stochastic integral $\int Z^i(\mbs)  \,dB^i(\mbs)$ is
defined in the sense of \cite{walsh}, with the $N$th coordinate
playing the role of the time variable and the other coordinates playing
the role of the spatial variables.

\begin{theorem}[(Cameron--Martin--Girsanov)]\label{thm3.3} If $(Z( \mbs),  \mbs\in
\IR_{+}^{N-1} \times[0, M])$ is such that $(L_{u}, u \in[0, M])$ is a
martingale with respect to the filtration $\ImG$, then the process
$(\tilde{B} (\mathbf{t}),  \mathbf{t} \in\IR_{+}^{N-1} \times[0,
M])$ defined by
\[
\tilde{B}(t_{1}, \ldots, t_{N}) = B(t_{1},
\ldots, t_{N}) - \int_{[0,t_1]\times\cdots\times[0, t_{N}]} Z(s_{1},
\ldots, s_{N}) \,ds_{1} \cdots \,ds_{N}
\]
is an $\IR^d$-valued Brownian sheet under the probability measure $Q$,
where $Q$ is defined by
\[
\frac{dQ}{dP} = L_{M}.
\]
\end{theorem}

We now fix $k \geq2$ and consider $k$ boxes $R_{1}, \ldots, R_{k}$ as
in the statement of Theorem~\ref{thm3.4}:
\[
R_{j} = \prod^{N}_{\ell=1}
\bigl[s^{0}_{j,\ell}, s^{1}_{j,\ell} \bigr],\qquad
j= 1, \ldots, k,
\]
where, for $\ell= 1, \ldots, N$, the intervals
\[
\bigl[s^{0}_{1,\ell}, s^{1}_{1,\ell}\bigr],
 \bigl[s^{0}_{2,\ell}, s^{1}_{2,\ell}\bigr]
, \ldots, \bigl[s^{0}_{k,\ell}, s^{1}_{k,\ell}
\bigr]
\]
are pairwise disjoint (i.e., the projection of the $R_{j}$ onto each
coordinate axis are pairwise disjoint). Without loss of generality, we
assume that
\[
s^{1}_{j-1, N} < s^{0}_{j, N},\qquad j=2,
\ldots, N
\]
(i.e., the projections of the $R_{j}$ onto the $N$th-coordinate axis are in increasing order).

Let
\begin{eqnarray*}
R &=& \Biggl( \prod^{N-1}_{\ell=1}
\bigl[s^{0}_{k,\ell}, s^{1}_{k,\ell}\bigr]
\Biggr) \times\bigl[s^{1}_{k-1,N}, s^{1}_{k,N}
\bigr],
\\
S& =& \Biggl( \prod^{N-1}_{\ell=1}
{\bigl]s^{0}_{k,\ell}, s^{1}_{k,\ell}\bigr[}^{c}
\Biggr) \times\bigl[0, s^{1}_{k-1, N}\bigr].
\end{eqnarray*}
Notice that $R_k \subset R$ and for $j=1,\ldots,k-1$, $R_j \subset S$.

\begin{lemma}\label{lem3.5} Let $M$ be as in Theorem~\ref{thm3.4}. There is a process
$(\hat{B}_{\mbt},  \mbt\in[0, M]^{N})$ with law mutually equivalent
to the law of $(B_{\mbt},  \mbt\in[0, M]^{N})$ such that
\[
\hat{B}(\mbt) = B(\mbt)\qquad \mbox{for } \mbt\in[0, M]^{N-1} \times
\bigl[0, s^{1}_{k-1, N}\bigr]
\]
and
\[
\hat{B}(\mbt) = B(\mbt)- E\bigl(B(\mbt) \mid\ImF(S)\bigr)\qquad\mbox{for } \mbt\in
R_{k}.
\]
In particular, $\hat{B}\vert_{R_{k}}$ and $ ( B\vert_{R_{1}}, \ldots
, B\vert_{R_{k-1}}  )$ are independent.
\end{lemma}

\begin{pf}
We apply Lemma~\ref{lem3.1} to the sets $R$ and $S$, yielding
the process $(\tilde{B}(\mbt),  \mbt\in R)$, such that $\tilde{B}(\mbt
) = E(B(\mbt) \mid\ImF(S))$, $\mbt\in R_{k}$. In particular, if we set
%
\begin{eqnarray}
\label{eq:13} \hat B(\mbt) &=& B(\mbt)  \qquad\mbox{for } \mbt\in{[0,
M]}^{N-1}\times\bigl[0, s^{1}_{k-1, N}\bigr],
\\
\hat B(\mbt) &= &B(\mbt) - \tilde B(\mbt) \qquad\mbox{for } \mbt\in
R_{k}, \label{eq:13a}
\end{eqnarray}
then $\hat B\vert_{R_{k}}$ and $(B\vert_{R_{1}}, \ldots,B\vert
_{R_{k-1}})$ are independent, since $B$ is a Gaussian process. The main
point of this lemma is to establish, after extending the definition of
$\hat B(\mbt)$ to $\mbt\in{[0, M]}^{N}$, that the law of $(\hat B(\mbt
),  \mbt\in{[0, M]}^{N})$ is mutually equivalent to the law of
$(B(\mbt),  \mbt\in{[0, M]}^{N})$.

For this, we will use Girsanov's theorem (Theorem~\ref{thm3.3}), by
constructing a process $(Z(\mbs))$ satisfying the assumption of Theorem~\ref{thm3.3} and such that
%
\begin{eqnarray}
\label{eq:13b} B(\mbt) - \int_{[0, t_{1}] \times\cdots\times[0, t_{N}]} Z(s_{1},
\ldots, s_{N})  \,ds_{1} \cdots \,ds_{N},
\nonumber
\\[-8pt]
\\[-8pt]
\eqntext{\mbt\in
\IR^{N-1} \times [0, M],}
\end{eqnarray}
agrees with $\hat B(\mbt)$ on $[0, M]^{N-1}\times[0, s^{1}_{k-1, N }]$
and on $R_{k}$. Using the formula in~\eqref{eq:13b} to define $\hat
B(\mbt)$ for all $\mbt\in\IR^{N-1} \times[0, M]$, this immediately
implies that the laws of $(\hat B(\mbt),  \mbt\in{[0, M]}^{N})$ and
$(B(\mbt),  \mbt\in{[0, M]}^{N})$ are mutually equivalent.

We note that for $\mbt= (t_{1}, \ldots, t_{N}) \in R$,
\[
\tilde B(\mbt)= \tilde B(t_{1}, \ldots, t_{N-1},
t_{N}) = \tilde B\bigl(t_{1}, \ldots, t_{N-1},
s^{1}_{k-1, N}\bigr),
\]
so $\tilde B(\mbt)$ does not depend explicitly on the $N$th-coordinate of $t$.

We now construct $Z(\mbs)$. Let
\[
U= \Biggl( \prod^{N-1}_{\ell=1} \bigl[0,
s^{1}_{k,\ell}\bigr] \Biggr) \times \bigl[s^{1}_{k-1, N},
s^{0}_{k, N}\bigr].
\]
We set
%
\begin{equation}
\label{eq:14a} Z(\mbs) \equiv0 \qquad\mbox{for } \mbs\notin U,
\end{equation}
and we define $Z(\mbs)$ for $\mbs\in U$ as follows. For $\mbt\in U
\cup R$, define
\[
p_{\ell} (\mbt) = s^{0}_{k,\ell} \vee
t_{\ell}, \qquad \ell= 1, \ldots, N-1,
\]
$p_{N}(\mbt) = s^{1}_{k-1, N}$, and
$
p(\mbt) = (p_{1}(\mbt), \ldots, p_{N}(\mbt))$.
Now let
%
\begin{equation}
\label{eq:14} 
\qquad F(\mbt)=\cases{ %
\displaystyle\frac{t_{N}- s^{1}_{k-1, N}}{s^{0}_{k, N} -s^{1}_{k-1, N}}
\Biggl( \prod^{N-1}_{\ell=1}
\frac{t_{\ell} \wedge s^{0}_{k,\ell}}{s^{0}_{k,\ell}}
\Biggr)\tilde B\bigl(p(\mbt)\bigr),&\quad $\mbox{if } \mbt\in U,$
\vspace*{2pt}\cr
0, &\quad $\mbox{otherwise,}$}
\end{equation}
so that $F(\mbt)$ is an $\IR^d$-valued multilinear interpolation of
$\tilde B(p(\mbt))$ with the process which vanishes on the coordinate
hyperplanes $1$ to $N-1$, and on the hyperplane $\IR^{N-1}\times\{
s^{1}_{k-1, N}\}$. In particular, for $\mbt\in U$,
%
\begin{equation}
\label{eq:15} F(\mbt) = 0\qquad \mbox{if } t_{1} =0 \mbox{ or } \cdots
\mbox { or } t_{N-1} =0 \mbox{ or } t_{N} =
s^{1}_{k-1, N}
\end{equation}
and
%
\begin{equation}
\label{eq:16} F(\mbt) = \frac{t_{N}- s^{1}_{k-1, N}}{s^{0}_{k, N} -s^{1}_{k-1, N}} \tilde B\bigl(t_{1},
\ldots, t_{N-1}, s^{1}_{k-1, N}\bigr)\qquad \mbox{if } \mbt
\in R.
\end{equation}
We note that $\mbt\mapsto F(\mbt)$ is piecewise $C^{\infty}$, and we set
\[
Z(s_{1}, \ldots, s_{N}) = \frac{\partial^{N}}{\partial{s_{1}}\cdots
\partial{s_{N}}}
F(s_{1}, \ldots, s_{N}).
\]
It is clear that $Z(\mbs)$ is a linear combination of the random
variables $B(s_{k,1}^{j (1)},\break \ldots,   s_{k,N-1}^{j(N-1)}, s^{1}_{k-1,
N})$ that come from Lemma~\ref{lem3.1}. Explicit formulas can be given,
for instance, letting $\dot{B}$ denote the white noise associated to $B$,
\begin{eqnarray*}
Z(\mbs) &=& \Biggl(\prod_{\ell=1}^{N-1}
\frac{1}{s^1_{k,\ell}-s^0_{k,\ell
}} \Biggr) \frac{1}{s^0_{k,N}-s^1_{k-1,N}}
\\
&&{} \times\dot{B} \bigl(\bigl[s^{0}_{k,1},
s^{1}_{k,1}\bigr] \times\cdots \times\bigl[s^{0}_{k, N-1},
s^{1}_{k, N-1}\bigr] \times\bigl[0, s^{1}_{k-1,
N}
\bigr]\bigr)\qquad \mbox{if } \mbs\in R,
\end{eqnarray*}
but we will not need them. We note, however, that $(Z(\mbs))$ is
adapted to the filtration $(\ImG_{u})$ defined in \eqref{fG}.

For $\mbt=(t_1,\ldots,t_N) \in\IR^{N}$, let
\[
\hat B(\mbt) = B(\mbt) - \int_{[0, t_{1}] \times\cdots\times[0,
t_{N}]} Z(s_{1},
\ldots, s_{N}) \,ds_{1} \cdots \,ds_{N}.
\]
Then \eqref{eq:13} is clearly satisfied by~\eqref{eq:14a}, and \eqref
{eq:13a} is satisfied since for $\mbt\in R_{k}$, by~\eqref{eq:15} and
\eqref{eq:16},
\begin{eqnarray*}
&& \int_{[0, t_{1}] \times\cdots\times[0, t_{N}]} Z(s_{1}, \ldots, s_{N}) \,ds_{1} \cdots \,ds_{N}
\\
&&\qquad = \int_{0}^{t_{1}} \,ds_{1} \cdots\int
_{0}^{t_{N-1}}\,ds_{N-1} \int
_{s^{1}_{k-1, N}}^{s^{0}_{k, N}}\,ds_{N} \frac{\partial^{N}}{\partial
{s_{1}} \cdots\partial{s_{N}}}
F(s_{1}, \ldots, s_{N})
\\
&&\qquad = \frac{s^{0}_{k,N}- s^{1}_{k-1,
N}}{s^{0}_{k,N}-s^{1}_{k-1,N}}\tilde B\bigl(t_{1}, \ldots, t_{N-1},
s^{1}_{k-1, N}\bigr)
\\
&&\qquad = \tilde B\bigl(p(\mbt)\bigr)
\\
&&\qquad = \tilde B(\mbt).
\end{eqnarray*}

In order to complete the proof, it remains to check that the
assumption of Theorem~\ref{thm3.3} is satisfied, and, in particular,
that the process
\[
L_{u}= \exp \biggl[\int_{\IR_{+}^{N-1} \times[0, u]} Z(\mbs)\cdot \,dB(
\mbs) - \frac{1}{2} \int_{\IR_{+}^{N-1} \times[0, u]} \bigl\Vert Z(\mbs )
\bigr\Vert^{2} \,d\mbs \biggr], \qquad u \in[0,M],
\]
is a martingale. Since $Z$ vanishes on $\IR^{N} \setminus U$, it
suffices, according to the extension of Novikov's criterion presented
in \cite{KS}, Chapter~3.5, Corollary~5.14, to check that for $n$
sufficiently large and $t_{i} = s^{1}_{k-1, N} + \frac{i}{n} (s^{0}_{k,
N}- s^{1}_{k-1, N})$, $i= 0, \ldots, n$,
\[
E \biggl[ \exp \biggl( \frac{1}{2} \int_{0}^{s^{1}_{k, 1}}
\,ds_{1} \cdots \int_{0}^{s^{1}_{k, N-1}}
\,ds_{N-1} \int_{t_{i-1}}^{t_{i}}\,ds_N \bigl\Vert Z(\mbs)\bigr\Vert^{2} \biggr) \biggr] < + \infty.
\]
But this follows from the fact that the integral is bounded by
\[
\frac{C}{n} \sup_{j \in\mathcal{J}_{N}}\bigl \Vert\bigl(B
\bigl(s_{k,1}^{j(1)}, \ldots, s_{k,N-1}^{j(N-1)},
s^{0}_{k-1,N}\bigr)\bigr)\bigr\Vert^{2}
\]
for some constant $C$ that depends only on $R_{k-1}$ and $R_{k}$, and
this random variable has a finite exponential moment if $n$ is
sufficiently large. The proof of Lemma~\ref{lem3.5} is complete.
\end{pf}

\begin{pf*}{Proof of Theorem~\ref{thm3.4}} We proceed by induction
on $k$. For $k=1$, there is nothing to prove. So, assume that $k \geq
2$ and that we have proved the statement for $k-1$.

We consider the two independent Brownian sheets $B$ and $W_{k}$. We
apply Lemma~\ref{lem3.5} to both of these processes, producing
processes $\hat B$ and $\hat{W}_{k}$ such that, in particular:
\begin{longlist}[(1)]
\item[(1)] $\hat B\vert_{R_{1}} = B\vert_{R_{1}}, \ldots,\hat B\vert
_{R_{k-1}} = B\vert_{R_{k-1}}$;
\item[(2)] $\hat B\vert_{R_{k}}$ and $(B\vert_{R_{1}}, \ldots, B\vert
_{R_{k-1}})$ are independent;
\item[(3)] $B\vert_{[0, M]^{N}}$ and $\hat B\vert_{[0, M]^{N}}$ have
mutually equivalent probability laws;
\item[(4)] $\hat{W}_{k}\vert_{R_{k}}$ and $W_{k}\vert_{R_{k}}$ have
mutually equivalent probability laws;
\item[(5)] $\hat B\vert_{R_{k}}$ and $\hat{W}_{k}\vert_{R_{k}}$ have
the \emph{same} probability law.
\end{longlist}
We write $\mathcal{L}(B\vert_{R_{1}}, \ldots, B\vert_{R_{k}})$ for the
probability law of the random vector $(B\vert_{R_{1}}, \ldots, B\vert
_{R_{k}})$, and use ``$\sim$'' to indicate mutually equivalent
probability laws. Then, by (3) and (1),
\begin{eqnarray*}
\mathcal{L}(B\vert_{R_{1}}, \ldots, B\vert_{R_{k}}) &\sim&
\mathcal{L} (\hat B\vert_{R_{1}}, \ldots, \hat B\vert_{R_{k-1}},
\hat B\vert _{R_{k}})
\\
& =& \mathcal{L}(B\vert_{R_{1}}, \ldots,B\vert_{R_{k-1}}, \hat B
\vert_{R_{k}}).
\end{eqnarray*}
By (2) and (5), and since $B$ and $W_k$ are independent,
\[
\mathcal{L}(B\vert_{R_{1}}, \ldots,B\vert_{R_{k-1}}, \hat B
\vert_{R_{k}}) = \mathcal{L}(B\vert_{R_{1}}, \ldots,B
\vert_{R_{k-1}}, \hat W_k\vert_{R_{k}}).
\]
Let $W_1,\ldots,W_{k-1}$ be independent Brownian sheets independent of
$W_k$ and $B$. Since $B$ and $W_k$ are independent, we can use the
induction hypothesis to see that
\[
\mathcal{L}(B\vert_{R_{1}}, \ldots,B\vert_{R_{k-1}}, \hat W
\vert_{R_{k}}) \sim\mathcal{L} (W_{1}\vert_{R_{1}},
\ldots, W_{k-1}\vert_{R_{k-1}}, \hat W_{k}
\vert_{R_{k}}).
\]
By (4) and the independence of $(W_1,\ldots,W_{k-1})$ and $W_k$, we
conclude that
\[
\mathcal{L} (W_{1}\vert_{R_{1}}, \ldots, W_{k-1}
\vert_{R_{k-1}}, \hat W_{k}\vert_{R_{k}}) \sim\mathcal{L}
(W_{1}\vert_{R_{1}}, \ldots, W_{k-1}
\vert_{R_{k-1}}, W_{k}\vert_{R_{k}}),
\]
%
and this proves Theorem~\ref{thm3.4}.
\end{pf*}


%

%




\printaddresses

\end{document}